\documentclass[11pt,reqno,palentino]{amsart} 

\usepackage{mypack}

\usepackage{latexsym}
\usepackage{a4wide}
\usepackage{amscd}
\usepackage{graphics}
\usepackage{amsmath}
\usepackage{amssymb}
\usepackage{mathrsfs}
\usepackage{accents}
\usepackage{enumerate}
\usepackage{url}
\usepackage[backgroundcolor=white,linecolor=gray]{todonotes}

\author{Thomas O. Rot}
\begin{document}
\title[Orientations]{The Morse-Bott inequalities, orientations, and the Thom isomorphism in Morse homology.}
\maketitle
\begin{abstract}
The Morse-Bott inequalities relate the topology of a closed manifold to the topology of the critical point set of a Morse-Bott function defined on it. The Morse-Bott inequalities are sometimes stated under incorrect orientation assumptions. We show that these assumptions are insufficient with an explicit counterexample and clarify the origin of the mistake. 
\end{abstract}

The Morse-Bott inequalities relate the topology of a closed manifold $M$ to the topology of the critical manifolds of a Morse-Bott function $f$ defined on it. Let $P_t(M)=\sum_i \mathrm{rank} \,H_i(M;\mZ) t^i$ denote the Poincar\'e polynomial with $\mZ$ coefficients, and let $MB_t(f)=\sum_j P_t(M_j)t^{| M_j|}$ be the Morse-Bott polynomial. Here the sum runs over the critical submanifolds of the function $f$, and $|M_j|$ denotes the index of the critical submanifold, i.e.\ the dimension of the negative normal bundle of $M_j$. We always take critical submanifolds to be connected components of the critical point set, which implies that the index is well defined. The Morse-Bott inequalities state that, under suitable orientation assumptions, there exists a polynomial $Q_t$ with non-negative coefficients such that 
\begin{equation}
\label{eq:morsebott}
MB_t(f)=P_t(M)+(1+t)Q_t.
\end{equation}
The orientation assumptions differ from paper to paper: In~\cite{jiang} no orientation assumptions are made. In~\cite{banyaga} it is assumed that the critical submanifolds are orientable and in \cite{banyagadynamical,hurtubise} it is additionally assumed that the ambient manifold is orientable. Below we show that these hypotheses are insufficient through an explicit counterexample, see also the errata~\cite{Banyaga:QUC7Szzy,Hurtubise:sz1P0qDm}. We construct a Morse-Bott function with orientable critical submanifolds on an orientable manifold that does not satisfy the Morse-Bott inequalities\ \bref{eq:morsebott}. 

Correct orientation assumptions for the Morse-Bott inequalities to hold are well known: Either one defines the Morse-Bott polynomial using homology with local coefficients in the negative normal bundles of the critical submanifolds, cf.\ \cite{bottmorsebott,Bismut:1986vx}, or one requires that the negative normal bundles are orientable as is for example done in\ \cite{Nicolaescu:2011jl}.

We now discuss the counterexample. Let $f\colon\thinspace\mathbb{R}P^5\rightarrow \mR$ be the function 
$$
f([x_0,x_1,\ldots,x_5])=\frac{-x_4^2+x_5^2}{x_0^2+\ldots+x_5^2}.
$$ 
The function is Morse-Bott with one minimum at $[0,0,0,0,1,0]$, one maximum at $[0,0,0,0,0,1]$, and one other critical submanifold $\mathbb{R}P^3=\{[x_0,x_1,x_2,x_3,0,0]\}$ of index $1$. We know that $P_t(\mR P^{2n+1})=1+t^{2n+1}$ and the Morse-Bott polynomial is
$$
MB_t(f)=P_t(\mathrm{min})t^0+P_t(\mR P^3)t^1+P_t(\mathrm{max})t^5=1+t+t^4+t^5.
$$
There exists no polynomial $Q_t$ with non-negative coefficients such that $1+t+t^4+t^5=1+t^5+(1+t) Q_t$, hence the Morse-Bott inequalities are violated. Note that odd dimensional projective spaces are orientable and that the normal bundle of $\mR P^3$ in $\mR P^5$ is also orientable. This example satisfies the hypotheses of~\cite{banyaga,banyagadynamical,hurtubise, jiang}, but violates Equation\ \bref{eq:morsebott}. 

The origin of the mistake is that the Thom isomorphism theorem requires homology with local coefficients or an orientation assumption on the bundle. We explain this briefly and refer to the appendix of~\cite{abbondandoloschwarz} for a more thorough discussion of the Thom isomorphism in Morse homology. Given a closed manifold $M$, a Morse function $f\colon\thinspace M\rightarrow \mR$, and vector bundles\footnote{Here and below we implicitly put suitable Riemannian metrics on the various bundles which are not really relevant to the discussion.}  $p^\pm\colon\thinspace E^\pm\rightarrow M$, define a Morse function $F\colon\thinspace E^+\oplus E^-\rightarrow \mR$ by 
\begin{equation}
\label{eq:morsefunction}
F(x^++x^-)=f(p^+(x^+))+\norm{x^+}_{E^+}^2-\norm{x^-}_{E^-}^2. 
\end{equation}
The critical points and the unstable manifolds of $F$ and $f$ are related as follows: $x\in \mathrm{crit}\, F$ if and only if $p(x)\in \mathrm{crit}\, f$ and $T_xW^u(x;F)\cong T_{p(x)}W^u(x;f)\oplus E^-_{p(x)}$. The generators of the Morse complexes of $f$ and $F$ coincide, but the grading is shifted by the dimension of $E^-$. To define the differential in Morse homology one orients the unstable manifolds. If the fibers $E^-_{p(x)}$ are coherently oriented, i.e.\ $E^-$ is oriented, then the manifolds $W^u(x;F)\cap W^s(y;F)$ and $W^u(p(x);f)\cap W^s(p(y);f)$ have the same induced orientation. The differentials of both Morse complexes then agree up to the shift in grading, and the Morse homology of $F$ is the Morse homology of $f$ shifted by the dimension of $E^-$. In general the Morse homology of $F$ computes the homology of the disc bundle in $E^-$ modulo its boundary, the sphere bundle in $E^-$, i.e.\ $HM_*(F)\cong H_*(DE^-,SE^-;\mZ)$. If the bundle $E^-$ is not orientable this does not relate to the Morse homology of $f$ directly. This sign issue can be explicitly seen for the standard Morse function on $M=S^1$ with two critical points and $E^-$ the non-orientable bundle over $S^1$. 

The failure of the Thom isomorphism to hold for non-orientable bundles comes up in Morse-Bott theory. To relate the homology of the critical submanifolds to the homology of the ambient manifold, the Morse-Bott function $f$ is perturbed to a Morse function $F=f+g$, where $g$ is sufficiently small and supported near the critical submanifolds. In a suitable tubular neighbourhood of a critical manifold the function $F$ has the form as in Equation\ \bref{eq:morsefunction}. Assuming that all the critical submanifolds have different values, the Morse complex of $F$ is filtered by choices of regular values of $f$. This gives a spectral sequence converging to the Morse homology of $F$. The second page of the spectral sequence are the groups isomorphic to $E_{i,j}^2\cong H_i(DN^-M_j,SN^-M_j;\mZ)$. If the negative normal bundles are orientable, this homology is the sum of the homologies of $M_j$'s shifted by the dimensions of the negative normal bundles $N^-M_j$'s. Standard commutative algebra gives the Morse-Bott inequalities.

\emph{Acknowledgement: The counterexample above is contained in the author's thesis~\cite{Rot:ww} written under supervision of Federica Pasquotto and Rob Vandervorst. The author thanks them, Hansj\"org Geiges, David Hurtubise, and Silvia Sabatini for correspondence on the topic of this paper.}

\bibliographystyle{abbrv} 
\bibliography{morse-bott}
\end{document}